\documentclass{amsart}  
\usepackage{amssymb}
\usepackage{graphicx}

\newtheorem{theorem}{Theorem}[section]

\newtheorem{problem}[theorem]{Problem}
\newtheorem{conjecture}[theorem]{Conjecture}
\newtheorem*{acknowledgement}{Acknowledgements}
\newtheorem*{prize}{Prize}
\theoremstyle{definition}
\newtheorem{definition}[theorem]{Definition}

\theoremstyle{remark}

\numberwithin{equation}{section}

\newcommand{\divrg}{\textrm{div}\,}

\title{A small collection of open problems}

\author{Giovanni Alessandrini}
\thanks{Trieste, Italy, 55gioale@gmail.com}

\begin{document}
\maketitle
\begin{abstract}{This paper collects some problems that I have encountered during the years, have puzzled me and which, to the best of my knowledge, are still open. Most of them are well--known and have been first  stated by other authors. In  this sad season of lockdown, I modestly try to contribute to  scientific interaction at a distance. Therefore all comments and exchange of information are most welcome.}
\end{abstract}



\section*{Introduction}

Here I collect a number of open problems that I have struggled with, and which I believe maintain some interest. I hope this collection may stimulate younger minds.

A masterly model, which gave me the inspiration to set up this small collection, has been the monumental Scottish Book \cite{ma}. A peculiarity which I especially liked in the Scottish Book is that often the proposers of the problems offered prizes for their solutions. I shall imitate this habit here. I promise the prizes to the solvers, the means of delivery are to be agreed depending on the occasion.

\section{Unique continuation for the $p$--laplacian}\label{sec:plap}
\begin{flushright}
\textit{\small Dedicato a Gisella, unica compagna della mia vita,\\ vicina e paziente  sempre, anche quando ``la matematica non funziona''.
}\end{flushright}

\medskip

Given $1<p<\infty$ we consider the so--called $p$--Laplace equation
\begin{equation}\label{plap}
\divrg (|\nabla u|^{p-2}\nabla u)  = 0
\end{equation}
in a connected open set $\Omega \subset \mathbb R^n$, $n\ge 2$. The natural space in which  we may search for weak (variational) solutions is $W^{1,p}(\Omega)$. The problem that I want to address here is whether
 solutions to \eqref{plap} satisfy the unique continuation property.
However, this formulation is somewhat vague, because the unique continuation property can be expressed in many ways.
More precise formulations are the following ones.
\begin{problem}[Weak unique continuation]\label{wuc}
Suppose $u$ solves \eqref{plap} and vanishes on a nonempty open set $U \subsetneq \Omega$, does $u\equiv 0$ on all of $\Omega$?
\end{problem}
\begin{problem}[Strong unique continuation]\label{suc}
Suppose $u$ solves \eqref{plap} and vanishes of infinite order at an interior point $x_0$ of $\Omega$, does $u\equiv 0$ on all of $\Omega$?
\end{problem}
If $p=2$ then \eqref{plap} reduces to the classical Laplace's equation, the solutions are harmonic functions which are well-known to be real analytic, and hence the \emph{strong} unique continuation property holds. If $n=2$ the strong unique continuation property holds true for any $p>1$. This was proved by Bojarski and Iwaniec \cite{bi} and by myself \cite{aplap} by two different methods,  a gap in the proof of Bojarski and Iwaniec  was later filled by Manfredi \cite{ma}.

Thus the problem, in both formulations \ref{wuc} and \ref{suc}, remains open in the case $n\ge3, p\neq 2$.

I should also mention that this problem was first proposed to me by Gene Fabes, in 1981, in an even stronger form:

\begin{problem}[Mukenhoupt]\label{muc}
Suppose $u$ solves \eqref{plap} and is not identically constant, is $|\nabla u|$ a Mukenhoupt weight?\end{problem}

For a definition of Mukenhoupt weights, and the basics of the theory, I refer to Coifman and Feffermann \cite{cf}.

Again, the answer is affirmative if $p=2$, this can be viewed within the general theory of unique continuation for linear elliptic equations, a benchmark in this respect is due to Nico Garofalo and F.H. Lin \cite{gl}. When $n=2$, the positive answer can be found in a joint paper with Daniela Lupo and Edi Rosset \cite{alr}.

\begin{prize} One bottle of barriqued Friuli Grappa.
\end{prize}



\section{Unique continuation along level surfaces}\label{sec:level}
\begin{flushright}
\textit{\small Dedicato a Giuseppe e Francesco Gheradelli, nonno e zio,\\ il ramo geometrico nel mio albero genealogico.
}\end{flushright}

This is probably the only problem of this collection for which I can claim some form of paternity. It originates from my work with Emmanuele Di Benedetto on the inverse problems of cracks \cite{adb}, and it could be used to extend some of those results to the case of cracks in inhomogeneous media. It was first formulated in a paper with Alberto Favaron \cite{afav}.

Consider an elliptic equation
\begin{equation}\label{ell}
\divrg (A\nabla u)  = 0
\end{equation}
in a connected open set $\Omega \subset \mathbb R^n$, $n\ge 2$. Here $A=A(x)$ is a symmetric matrix valued function, satisfying uniform ellipticity
\begin{equation}\label{ellip}
K^{-1}|\xi|^2\le A(x)\xi \cdot \xi \le K |\xi|^2 \ , \forall x, \xi \in \mathbb R^n \ ,
\end{equation}
 and Lipschitz continuity
 \begin{equation}\label{lip}
|A(x)-A(y)|\le E|x-y| \ , \forall x, y \in \mathbb R^n \ , 
\end{equation}
for some positive constants $K, E$. These  assumptions are stated in order to guarantee that the standard unique continuation property applies, see again \cite{gl}. I propose the following less standard form of unique continuation.
\begin{problem}[Unique continuation along level surfaces]\label{luc}
Let $u,v$ be two noncostant solutions to \eqref{ell}, let $S\subset \Omega$ be a connected $(n-1)$--dimensional smooth hypersurface, and let $\Sigma$  be an open  proper subset of $S$. Suppose $u=0$ on $S$ and $v=0$ on $\Sigma$. Is it true that $v=0$ on all of $S$?\end{problem}

If $A\equiv I$, then $u, v$ are harmonic and in this case, $S$ is an analytic hypersurface, see for instance \cite[Appendix A]{adb}, hence the answer is affirmative, by analytic continuation. Some smoothness on $S$ must be indeed assumed. This can be easily seen in the harmonic two--dimensional case:

Let $\Omega = \mathbb R^2$, $u=xy$, $v=y$. In this case, we may pick  $S= \{x=0, y\ge 0\}\cup \{x\ge0, y=0\}$, which is a simple curve with a corner point at the origin and $\Sigma=\{x>0, y=0\} \subset S$ is the positive horizontal semiaxis. Clearly $v$ vanishes on $\Sigma$, but it does not vanish on all of $S\subset \{u=0\}$.

From another point of view, it may be noted that a parallel could be drawn with the issue of unique continuation at fixed time   for parabolic equations, which was treated jointly with Sergio Vessella \cite{ave}.

\begin{prize} Three Havana cigars.
\end{prize}

\section{The inclusion problem}\label{sec:incl}

Consider a connected open set $\Omega \subset \mathbb R^n$, $n\ge 2$, with smooth boundary: Let $k>0, k\neq 1$ be given, and let $D$ be an open set compactly contained in $\Omega$. It is well--known that the Dirichlet problem
\begin{equation}\label{incl}
\left\{
\begin{array}{lll}
\divrg((1 + (k-1) \chi_D) \nabla u) =0,&\hbox{in}&\Omega,\\
u=\varphi,&\hbox{on}&\partial \Omega.
\end{array}
\right.
\end{equation}
This problem models the distribution of an electrostatic potential $u$ in a body $\Omega$ with homogeneous conductivity $1$  which contains in its interior   an inclusion $D$ whose conductivity is $k\neq 1$. 

Such a problem is well--posed in the Sobolev space $W^{1,2}(\Omega)$ for any given Dirichlet data $\varphi$ in the trace space $W^{\frac{1}{2},2}(\partial\Omega)$. As a consequence, the so--called Dirichlet--to--Neumann map
\[\Lambda_D: W^{\frac{1}{2},2}(\partial\Omega) \ni \varphi \to \nabla u \cdot \nu \in W^{-\frac{1}{2},2}(\partial\Omega)  \]
is well defined, here $\nu$ denotes the outer unit normal to $\partial\Omega$.

A special case of the celebrated Calder\'on's inverse conductivity problem \cite{cal} is to determine $D$, given $\Lambda_D$. Victor Isakov proved uniqueness \cite{is88}.

Loosely speaking, the boundary of $D$ is determined by $n-1$ parameters, whereas the Dirichlet--to--Neumann map can be associated to a function depending on $2(n-1)$ variables. Therefore the problem of determining $D$, from $\Lambda_D$ looks dimensionally overdetermined. It is thus interesting to see if a limited sample of $\Lambda_D$  might suffice to uniquely determine $D$. A review of the available results can be found in \cite{acat}.

Let us recall here in particular a result by Jin Keun Seo  in dimension $n=2$ \cite{seo}. Assuming a priori that $D$ is a simply connected polygon, given two \emph{cleverly} chosen Dirichlet data $\varphi_1, \varphi_2$, then the corresponding Neumann data $\nabla u_1\cdot \nu , \nabla u_2\cdot \nu$ uniquely determine $D$. Here $u_i$ denotes the solution to \eqref{incl} when $\varphi=\varphi_i$, $i=1,2$.

What is a clever choice of $\varphi_1, \varphi_2$? 

\begin{definition} \label{clev}
We say that the pair of Dirichlet data $\varphi_1, \varphi_2$ satisfies Rad\'o's condition if the mapping $\Phi = (\varphi_1, \varphi_2): \partial \Omega \to \Gamma$ is a homeomorphism onto the boundary $\Gamma$ of an open convex set $G\subset \mathbb R^2$.
\end{definition}
 Seo's uniqueness theorem was originally stated under slightly different conditions, but it can be readily seen that his arguments apply equally well under the condition just stated. The rationale behind this \emph{clever choice} of Dirichlet data, is that under such condition, the mapping $U=(u_1, u_2)$ becomes a homeomorphism of $\Omega$ onto $G$, see my joint paper with Enzo Nesi \cite{an}. In fact such a condition was first devised by Rad\'o in the context of planar harmonic mappings, see for instance Duren \cite{du}.
 
 \begin{problem}[The inclusion problem with finite data]\label{pro:inc}
Let $n=2$, let us assume a priori  $D$ be  simply connected and with smooth boundary. Let the pair of Dirichlet data $\varphi_1, \varphi_2$ satisfy  Rad\'o's condition, is $D$ uniquely determined by the Neumann data $\nabla u_1\cdot \nu , \nabla u_2\cdot \nu$?
\end{problem}
 One clue that Rad\'o's condition might be the appropriate one comes from the fact that it has shown to be effective in the germane inverse problem of cracks in dimension 2. The main results, and more references, can be found in the paper by Alvaro Diaz Valenzuela and myself \cite{ad}, and in the one by Jin Keun Seo with Kim \cite{ks}.
 
 Of course it would be meaningful to formulate Problem \ref{pro:inc} also in higher dimensions. However it is not at all clear what might be a  clever choice of Dirichlet data in such a case.
 
 \begin{prize} Three bottles of Friulano (once known as Tocai).
\end{prize}

\section{The crack problem}\label{sec:crack}

\begin{flushright}
\textit{\small Dedicato a Damien, che mi d\`a la forza di sperare nel futuro. 
}\end{flushright}

\medskip

Consider, in a bounded and connected open set $\Omega \subset \mathbb R^3$ a two--dimensional orientable smooth surface $\Sigma \Subset \Omega$ having a  simple closed curve as  boundary. We view $\Sigma$ as a fracture in a homogeneous electrically conducting body $\Omega$. If we prescribe a stationary current density $\psi$ (having zero average) on $\partial \Omega$ then the electrostatic potential $u$ is governed by the following boundary value problem (suitably interpreted in a weak sense)
\begin{equation}\label{insul}
\left\{
\begin{array}{lll}
\Delta u =0,&\hbox{in}&\Omega\setminus \Sigma,\\
\nabla u\cdot \nu=0 , &\hbox{on either side of}& \Sigma, \\
\nabla u\cdot \nu=\psi,&\hbox{on}&\partial \Omega.
\end{array}
\right.
\end{equation}
The inverse problem is:
\begin{problem}[The insulating crack problem]\label{pro:cra}
To find appropriately chosen current profiles $\psi_1, \ldots \psi_N$ such that, letting $u_1, \ldots u_N$ be the corresponding potentials, the boundary measurements $u_1|_{\partial \Omega}, \ldots u_N|_{\partial \Omega}$ uniquely determine $\Sigma$.
\end{problem}
When $\Sigma$ is a portion of a plane, Kubo \cite{ku} found a triple of suitable current profiles for which uniqueness holds, 
in \cite{adb} suitable pairs were found. 

When $\Sigma$ is allowed to be curved and the full Neumann--to--Dirichlet map
\[ \mathcal N_{\Sigma}: \psi \to u|_{\partial \Omega}\]
is known, uniqueness was proven by Eller \cite{ell}.

\begin{prize} One bottle of white Friuli Grappa.
\end{prize}

\section{Payne's nodal line conjecture}\label{sec:payne}

In this Section and in the following one I shall discuss problems arising from Courant's Nodal Domain Theorem \cite{ch}. A brief introduction may be useful. Consider, in a bounded and connected open set $\Omega \subset \mathbb R^n$, $n\ge 2$, the eigenvalue problem
\begin{equation}\label{eigen}
\left\{
\begin{array}{lll}
\Delta u + \lambda u=0,&\hbox{in}&\Omega,\\
u=0,&\hbox{on}&\partial \Omega.
\end{array}
\right.
\end{equation}
It is well known that there exists a complete set of eigenfunctions $\{u_m\}$, orthonormal in $L^2(\Omega)$, with corresponding eigenvalues $\{\lambda_m\}$. The eigenvalues are all positive and can be arranged in a nondecreasing, diverging sequence $0< \lambda_1 < \lambda_2\le \lambda_3\le\ldots.$.

Courant's Nodal Domain Theorem states that every eigenfunction $u_m$ corresponding to the $m$--th eigenvalue $\lambda_m$ has at most $m$ nodal domains, that is the set $\{x\in \Omega| u_m(x) \neq 0\}$ has at most $m$ connected components.

In particular, $u_1$ has constant sign, and any second eigenfunction $u_2$ (which, being orthogonal to $u_1$, must change its sign) has exactly two nodal domains.

L. E. Payne \cite[Conjecture 5]{pay} posed the following conjecture.

\begin{conjecture}[Payne]\label{payne}
Let $n=2$, for any second eigenfunction $u_2$ the nodal line \[\overline{\{x\in {\Omega} | u_2(x) = 0\}}\] is not a closed curve.
\end{conjecture}

Melas \cite{me} proved the conjecture if $\Omega$ is convex and has smooth boundary, the smoothness assumption was removed in \cite{ga}. Hoffmann-Ostenhof,  Hoffmann-Ostenhof and Nadirashvili \cite{hona} showed by an example that the conjecture may fail if $\Omega$ is not simply connected.

The following remains an open problem.

 \begin{problem}\label{pro:payne}
Assume  $\Omega$ be  simply connected and with smooth boundary. Prove that the nodal line of any second eigenfunction is a simple open curve whose endpoints are two distinct points of $\partial \Omega$.
\end{problem}

It must be noted that the conjecture also fails if the eigenvalue problem \eqref{eigen} is slightly modified by replacing the Laplacian $\Delta$ with an operator of the form $\Delta -q$ where $q$ is a variable coefficient, see the example by Lin and Ni \cite{linni}.

\begin{prize} Three bottles of Terrano.
\end{prize}

\section{Extending Courant's nodal domain thorem}\label{sec:courant}
Let us consider again an elliptic eigenvalue problem, but now we admit variable coefficients:
\begin{equation}\label{eigenvar}
\left\{
\begin{array}{lll}
\divrg (A\nabla u)  - qu+ \lambda \rho u=0,&\hbox{in}&\Omega,\\
u=0,&\hbox{on}&\partial \Omega.
\end{array}
\right.
\end{equation}
We shall assume throughout the ellipticity condition \eqref{ellip}, and also
\[ q, \rho \in L^{\infty}(\Omega)\ ,  \rho\ge K^{-1}> 0 \ . \]
It is known that Courant's Nodal Domain  Theorem maintains its validity under one of the following conditions, either
$n=2$ or $A$ satisfies the Lipschitz condition \eqref{lip}. See \cite{anod} for proofs and bibliography. In fact under such assumptions, solutions to
\begin{equation}\label{elleq}
\divrg (A\nabla u)  - qu+ \lambda \rho u=0
\end{equation}
satisfy the unique continuation property. This is in fact a typical ingredient in the proof of the Nodal Domain Theorem.

If $n\ge3$ and \eqref{lip} is relaxed to a H\"older condition
 \begin{equation}\label{lhold}
|A(x)-A(y)|\le E|x-y| ^{\alpha}\ , \forall x, y \in \mathbb R^n \ , 
\end{equation}
for some $\alpha \in (0,1)$ and some positive constant $E$, then a weaker result is known: any $m$--th eigenfunction $u_m$ has at most $2(m-1)$ nodal domains, for every $m\ge 2$. See \cite[Theorem 4.5]{anod}.
Recall that if $A$ is merely H\"older, then unique continuation may fail \cite{pli, mi1, mi2}.
 \begin{problem}[Courant]\label{pro:courant}
Let $n \ge 3$. Under which conditions on $A$ (not implying the unique continuation property for \eqref{elleq}) Courant's Nodal Domain Theorem remains valid?
\end{problem}

\begin{prize} Three bottles of Vitovska.
\end{prize}


\section{The final problem}\label{sec:final}

\begin{flushright}
\textit{\small Dedicato
a Luigi Gherardelli, lo zio Gigi, ingegnere idraulico. 
}\end{flushright}

\medskip

Consider, for a bounded, connected open set $\Omega \subset \mathbb R^n$, the following initial--boundary value problem for a parabolic equation
\begin{equation}\label{parab}
\left\{
\begin{array}{lll}
u_t- \divrg(a \nabla u) =f,&\hbox{in}&\Omega\times (0,T),\\
u=0,&\hbox{on}&\partial \Omega\times (0,T) \cup \Omega\times \{0\}.
\end{array}
\right.
\end{equation}
Here $a=a(x)$ is assumed to be a time independent scalar coefficient, satisfying uniform ellipticity
\[ K^{-1}\le a(x)\le K  \ .\]
In \cite{avess} Sergio Vessella and I formulated the following inverse problem.
\begin{problem}[Parabolic inverse problem]\label{pro:final1}
Let $n=2$. Given $f$ on all of $\Omega\times (0,T)$, and given $u(x,t)$ for all $x\in \Omega$ and for one or more fixed values of $t>0$, is $a$ uniquely determined?\end{problem}
The limitation of the space dimension to $n=2$ was motivated by applications, more specifically: identification of transmissivity $a$ in groundwater flow by piezometric head measurements $u$. The origin of the problem should be searched in the hydrogeology literature, and it might be dated much earlier, see for instance \cite{yd}.

The problem however maintains its interest in any dimension, and is still open in general. A solution under some special assumptions  on boundary data and when $n=1$ was given by Victor Isakov in \cite{is91}.

Here I propose a variant to Problem \ref{pro:final1}.
\begin{problem}[Parabolic inverse problem with final data]\label{pro:final2}
Consider the following initial--boundary value parabolic problem
\begin{equation}\label{parabh}
\left\{
\begin{array}{lll}
u_t- {\rm div}(a \nabla u) =0,&\hbox{in}&\Omega\times (0,T),\\
u=g(x,t),&\hbox{on}&\partial \Omega\times (0,T) , \\
u=h(x) ,&\hbox{on}&\Omega\times \{0\}.
\end{array}
\right.
\end{equation}
Given appropriately chosen functions $g$ and $h$, and assuming $a$ known on $\partial \Omega$, does the knowledge of $u(\cdot, T)$ uniquely determines $a$?
\end{problem}
This problem can be seen as an extension to a non--stationary setting of an inverse elliptic problem with interior data, for which uniqueness can be proven \cite{a86, a2016}.  In such an elliptic context, it can be seen that some boundary information on the coefficient $a$ is needed, that is why I have added that kind of information also in the parabolic question.

\begin{prize} Three bottles of Ribolla.
\end{prize}

\begin{acknowledgement}
This note is intended as a token of gratitude for my colleagues of the Dipartimento di Matematica e Geoscienze with whom I have shared almost three decades of work. 
I wish 
also 
to thank all the friends and colleagues with whom I had the privilege to collaborate in the course of my career.
A special thank goes to my friends Rolando Magnanini, Enzo Nesi, Sergio Vessella, whose advice on a first draft of this paper I treasured.
\end{acknowledgement}

 \end{document}